\documentclass[twoside,12pt,a4paper,mathscr]{article}
\usepackage{amsmath,amssymb,amsthm,euscript,amscd,a4wide}
\usepackage[all]{xypic}
\CompileMatrices
\newtheorem{thm}{Theorem}[section]
\newtheorem*{teor}{Theorem}

\newtheorem{lemma}[thm]{Lemma}
\newtheorem{prop}[thm]{Proposition}
\newtheorem{defin}[thm]{Definition}

\newtheorem{rem}[thm]{Remark}
\newtheorem{ex}[thm]{Example}

\newenvironment{remark}{\begin{rem}\rm}
       {\hfill$\vartriangle$\mbox{\hskip3pt}\end{rem}}
\let\cal=\mathcal
\newcommand{\Pc}{{\cal P}}

\newcommand{\Z}{{\mathbb Z}}

\newcommand{\lra}{\longrightarrow}

\newcommand{\Ext}{\operatorname{Ext}}
\newcommand{\Ker}{\operatorname{Ker}}
\newcommand{\Ima}{\operatorname{Im}}
\newcommand{\Coker}{\operatorname{Coker}}

\newcommand\ch{\operatorname{ch}}
\newcommand\WIT{\operatorname{WIT}}
\newcommand\IT{\operatorname{IT}}
\newcommand\rk{\operatorname{rk}}

\newcommand\Id{\operatorname{Id}}

\newcommand{\dual}{\vee}
\newcommand{\isom}{\cong}
\newcommand{\HH}{\mathscr H}
\newcommand{\mat}[4]{\left( \begin{array}{cc} #1 & #2 \\ #3 & #4 \end{array} \right)}
\newcommand{\col}[2]{\left( \begin{array}{c} #1 \\ #2 \end{array} \right)}
\let\tensor=\otimes

\newcommand{\QQ}{\mathcal Q}

\newcommand{\Ltensor}{\mathbin{\buildrel{\mathbf L}\over{\tensor}}}

\newcommand{\bysame}{$\raise.2em\hbox to 3em{\hrulefill}$\thinspace, }
\begin{document}
\thispagestyle{empty}
\begin{center}
{\bfseries\Large  Fourier-Mukai transforms for coherent systems
 on elliptic curves}
\par\addvspace{20pt}
{\sc Daniel Hern\'andez Ruip\'erez}
and {\sc Carlos Tejero Prieto}
\par\medskip
Departamento de Matem\'aticas and Instituto de F\'{\i}sica Fundamental y Matem\'aticas\\
Universidad de Salamanca\\ Plaza de la Merced 1-4, 37008
Salamanca, Spain
\end{center}
\vfill
\begin{quote} \footnotesize {\sc Abstract.}
We determine all the Fourier-Mukai transforms for coherent systems
consisting of a vector bundle over an elliptic curve and a
subspace of its global sections, showing that these transforms are
indexed by the positive integers. We prove that the natural
stability condition for coherent systems, which depends on a
parameter, is preserved by these transforms for small and large
values of the parameter. By means of the Fourier-Mukai transforms
we prove that certain moduli spaces of coherent systems
corresponding to small and large values of the parameter are
isomorphic. Using these results we draw some conclusions about the
possible birational type of the moduli spaces. We prove that for a
given degree $d$ of the vector bundle and a given dimension of the
subspace of its global sections there are at most $d$ different
possible birational types for the moduli spaces.
\end{quote}
\vfill \leftline{\hbox to8cm{\hrulefill}}\par {\footnotesize
\noindent  This research has been partially supported by the
research projects MTM2006-04779 of the Spanish MEC, SA114/04 of
the ``Junta de Castilla y Le\'on'' and by the European Scientific
Exchange Programme ``Geometry and topology of moduli spaces'' of
the Royal Society of London and the Spanish CSIC under grant
15646.
\\
\noindent\emph{E-Mail addresses:} {\tt  ruiperez@usal.es, carlost@usal.es} \\
\noindent\emph{Mathematics Subject Classification:}
14D20, 14H60,14J60 and 14H52. \\
\noindent To appear in \emph{Journal of the London Mathematical Society}}

\eject\thispagestyle{empty}\mbox{\ \ \ }
\setcounter{page}{1}

\section{Introduction}

Fourier-Mukai transforms have proved to be a very useful tool for
studying moduli spaces of vector bundles. These transforms have
been used in a wide range of situations that include, among
others, abelian varieties \cite{Mu1}, abelian and K3 surfaces
\cite{Mac, Yo, BBH2}, and elliptic fibrations \cite{BBHM,Br,HP}.
The applications to moduli spaces of augmented vector bundles has
been initiated recently and cover the cases of Higgs bundles
\cite{BB} and holomorphic triples \cite{GHPT}.

In this paper we endeavor to apply Fourier-Mukai transforms to
study the moduli spaces of a different type of augmented
structure, namely coherent systems on elliptic curves. Their
moduli spaces have been studied by Lange and Newstead in the
recent paper \cite{LN}. Although coherent systems may be thought
of as particular instances of holomorphic triples, their stability
notions are actually different. Therefore, the results of this
paper do not follow from the study carried out in \cite{GHPT}.
Moreover, the scope of the present work is wider than that of
\cite{GHPT} in two respects. Firstly, we deal with all the
Fourier-Mukai transforms that act naturally on the category of
coherent systems rather than concentrating on the original
transform introduced by Mukai. Secondly, we do not restrict
ourselves to any coprimality assumptions on the topological
invariants of the coherent systems.

Let us explain our method. As a first step, which is necessary in
order to study all the Fourier-Mukai transforms that act on
coherent systems, we must have a clear picture of ordinary
Fourier-Mukai transforms on an elliptic curve $X$ and their action
on bundles. Based on previous results of Bridgeland \cite{Br} and
Hille-van den Bergh \cite{Hi}, we obtain a complete
characterization of the action of any Fourier-Mukai transform on
stable and semistable bundles on an elliptic curve.

A coherent system $(E,V)$ of type $(r,d,k)$ consists of a vector
bundle $E$ with rank $r$, degree $d$ and a $k$-dimensional
subspace $V$ of its global sections. As is well known, $(E,V)$ can
be identified with its evaluation map $V\otimes\mathcal O_X\to E$,
and in this way we may regard coherent systems as a subcategory
$S(X)$ of the abelian category $\mathcal C(X)$ whose objects are
arbitrary maps $\varphi\colon V\otimes\mathcal O_X\to\mathcal E$
where $V$ is a finite dimensional vector space and $\mathcal E$ is
any coherent sheaf. The second step of our method aims at the
determination of all the Fourier-Mukai transforms that act on
coherent systems on an elliptic curve. Since any Fourier-Mukai
transform $\Phi$ acts in a natural way on arbitrary sheaf maps, we
only have to determine which transforms leave the category $S(X)$
invariant. We start by considering the same, but easier, problem
for $\mathcal C(X)$. Once this is solved, the spectral sequence
$E_2^{p,q}=H^p(X,\Phi^q(\mathcal{E}))\Rightarrow
H^{p+q}(Y,\Phi(\mathcal{E}))$ attached to any integral functor
$\Phi$ allows us to obtain the Fourier-Mukai transforms that act
on $S(X)$ (see Propositions \ref{fm-cs-I}, \ref{psi-it1}).

With these results at hand we proceed to study the preservation of
stability under Fourier-Mukai transforms. Let us recall that the
stability notion for coherent systems depends on a real parameter
$\alpha>0$. Moreover, the $\alpha$-range for which $\alpha$-stable
coherent systems do exist is divided into a finite number of open
intervals such that the moduli spaces $G(\alpha;r,d,k)$ for any two
values of $\alpha$ inside the same interval coincide. The behavior
of the bundle that forms a coherent system with respect to
Fourier-Mukai transforms is not always easy to determine. However,
there are two moduli spaces $G_0(r,d,k)$, and $G_L(r,d,k)$ if we
assume that $k<r$, (corresponding to the first and last open
intervals, respectively) for which we are able to establish the
preservation of stability. Our study is based on the fact that in
any coherent system $(E,V)\in G_0(r,d,k)$ over a projective curve
the bundle $E$ is semistable, whereas if $(E,V)\in G_L(r,d,k)$ then
$E$ is the middle term of a BGN extension (see Definition \ref{bgn})
with semistable quotient. The key point now is the determination of
the sufficient conditions for the converses to be true, see
Propositions \ref{esta}, \ref{esta-bgn}. In the particular case of
the moduli space $G_L(r,d,k)$, our techniques allow us to prove in
addition the preservation of the families $\mathcal{BGN}(r,d,k)$,
$\mathcal{BGN}^s(r,d,k)$ of BGN extension classes of type $(r,d,k)$
on $X$ in which the quotient is semistable or stable, respectively.
These families are important on their own and are related to the
moduli space by the chain of inclusions
$\mathcal{BGN}^s(r,d,k)\hookrightarrow G_L(r,d,k)\hookrightarrow
\mathcal{BGN}(r,d,k)$.

These ideas allow us to prove the main results of the paper
(Theorems \ref{sm-iso}, \ref{lg-iso}), namely that certain moduli
spaces of coherent systems are isomorphic, with the isomorphisms
given by the Fourier-Mukai transforms $\Phi_a$ described in
Section \ref{fm-coh-sys}. More precisely we have:

\begin{teor}
Let $a$ be any positive integer. The Fourier-Mukai transforms
$\Phi_a$ induce isomorphisms of moduli spaces

\begin{enumerate}
\item $\Phi_a^0\colon G_0(r,d,k)\to G_0(r+ad,d,k)$

\item $\Phi_a^0\colon G_L(r,d,k)\to G_L(r+ad,d,k)$ ($k<r$)
\end{enumerate}
\end{teor}

It follows that the isomorphism type of $G_0(r,d,k)$ and
$G_L(r,d,k)$ (with $k<r$) depends only on the class $[r]\in\mathbb
Z/d\,\mathbb Z$.

In the particular case where $r$ and $d$ are coprime, Lange and
Newstead proved (see \cite[Proposition 8.1]{LN}) that the moduli
space $G_0(r,d,k)$ can be identified with the Grassmannian bundle
over the moduli space of stable bundles $M(r,d)$ of subspaces of
dimension $k$ in the fibres of the Picard bundle $p_{2*}\mathcal
U$, where $\mathcal{U}$ is a Poincar\'{e} bundle on $X\times M(r,d)$
and $p_2$ is the projection onto the second factor. In the same
way if $r-k$ and $d$ are coprime with $0<k<r$, they showed that
$G_L(r,d,k)$ can be identified with the Grassmannian bundle over
$M(r-k,d)$ of subspaces of dimension $k$ in the fibres of the
Picard bundle $(p_{2*}{\mathcal U}^\prime)^*$, where
$\mathcal{U}^\prime$ is a Poincar\'{e} bundle on $X\times M(r-k,d)$.
Therefore, in these particular cases it is conceivable that one
could prove that $G_0(r,d,k)$ and $G_0(r+ad,d,k)$ are isomorphic,
as well as $G_L(r,d,k)$ and $G_L(r+ad,d,k)$, by solely using the
properties of the moduli spaces of stable bundles and of their
Poincar\'{e} bundles. However it is interesting to point out that our
proof gives explicit isomorphisms by means of the Fourier-Mukai
transforms.

As an application of our results we are able to draw some
conclusions about the possible birational type of the moduli spaces.
Lange and Newstead proved in \cite{LN} that the birational type of
$G(\alpha;r,d,k)$ is independent of $\alpha$. Therefore, the
birational type of the moduli spaces can be determined by
considering any of the finitely many different moduli spaces. Taking
into account the preceding Theorem, the natural candidates to
consider are $G_0(r,d,k)$ and $G_L(r,d,k)$. However, the latter is
less convenient than the former since we have studied it under the
restriction $k<r$. In this way we prove that the birational types of
$G(\alpha;r,d,k)$ and $G(\alpha;r+ad,d,k)$ are the same. Thus, the
birational type of $G(\alpha;r,d,k)$ depends only on the class
$[r]\in\mathbb Z/d\,\mathbb Z$. We conclude that for fixed $d$ and
$k$, there are at most $d$ different birational types for the moduli
spaces of coherent systems of type $(r,d,k)$.

We finish with a brief comment on the possible generalizations and
applications of the work carried out in this paper. Our method,
based on the Fourier-Mukai transform, is not suitable for other
than elliptic curves. Fourier-Mukai transforms have proved to be
well-behaved in the case of varieties with trivial canonical sheaf
(and some other closely related cases) and in this sense, one
natural generalization of our paper would be to the case of
coherent systems on higher dimensional abelian varieties. Other
possible future application is to relative coherent systems on
elliptic fibrations, following the guidelines established by many
authors (see for instance \cite{AHR,HP,BBH2}), in the same vein of
the Friedman-Morgan-Witten classical application of sheaves on
elliptic fibrations to mirror symmetry. We plan to study these new
cases in forthcoming papers.

The paper is organized as follows. In Section \ref{ftransform} after
briefly recalling some facts on integral functors we obtain a
complete description of Fourier-Mukai transforms on an elliptic
curve and their action on stable and semistable bundles. Section
\ref{trans-coh-sys} is devoted to the determination of all the
Fourier-Mukai transforms that act on coherent systems on an elliptic
curve. In Section \ref{pres-stab} we address the question of
preservation of stability of coherent systems under the
Fourier-Mukai transforms determined in Section \ref{trans-coh-sys}.
The study of this question is divided in two cases, corresponding to
the moduli spaces for small and large values of the parameter
$\alpha$, which require different techniques. We prove in this
section the isomorphisms induced by the Fourier-Mukai transforms
between moduli spaces, Theorems \ref{sm-iso}, \ref{lg-iso}. The
isomorphism between families of BGN extensions classes is proved in
Theorem \ref{iso-bgn}. Finally, in Section \ref{bir-type} we apply
the previous results to obtain some conditions for the birational
type of the moduli spaces of coherent systems.

In this paper we work over the field of complex numbers $\mathbb
C$.

\section{Fourier-Mukai transforms on elliptic cur\-ves}
\label{ftransform}

We begin this section by briefly recalling the facts on integral
functors that we use in the paper, for further details the reader is
referred, for instance, to \cite{Mu,BO,AHR,BBHJ}. Given two smooth
varieties $X$ and $Y$ and an object $\Pc\in D(X\times Y)$, one
defines a functor $\Phi^\Pc_{X\to Y}: D(X)\lra D(Y)$ between the
bounded derived categories of $X$ and $Y$, by the formula
$$\Phi^\Pc_{X\to Y}(-)=\mathbf{R}\pi_{Y,*}(\pi^*_X(-)\Ltensor\Pc),$$
where $\pi_X$ and $\pi_Y$ are the projections from $X\times Y$ to
$X$ and $Y$, respectively. We say that $\Phi=\Phi^\Pc_{X\to Y}$ is
an integral functor with kernel $\Pc$.

\begin{defin}
Given an object $E$ of $ D(X)$ we set $\Phi^i(E)=\HH^i(\Phi(E)).$
A sheaf $E$ on $X$ is said to be $\Phi$-$\WIT_i$ if $\Phi^j(E)=0$
for all $j\neq i$. We say $E$ is $\Phi$-$\WIT$ if it is
$\Phi$-$\WIT_i$ for some $i$, and in this case we often write
$\Hat E$ for $\Phi^i(E)$, and refer to $\Hat E$ as the {\it
transform} of $E$.
\end{defin}

Define
$$\QQ=\Pc^{\dual}\tensor\pi_Y^*\omega_Y,$$
and put $\Psi=\Phi^{\QQ}_{Y\to X}$. It is a simple consequence of
Grothendieck-Verdier duality that $\Psi[\dim Y]$ is a left adjoint
of $\Phi$. Moreover, if $\Phi$ is fully faithful one has an
isomorphism of functors
$$\Psi\circ\Phi\isom \Id_{ D(X)}[-\dim Y].$$

\begin{defin}
We say that the integral functor $\Phi\colon D(X)\to D(Y)$ is a
Fourier-Mukai functor if it defines an equivalence of categories. If
in addition, the kernel $\Pc$ of $\Phi$ is a sheaf, then we say that
$\Phi$ is a Fourier-Mukai transform.
\end{defin}

Now we recall the definition of $\IT$-sheaves for a Fourier-Mukai
transform.

\begin{defin}
Let $\Phi$ be a Fourier-Mukai transform. A sheaf $E$ on $X$ is said
to be $\Phi$-$\IT_i$ if it is $\Phi$-$\WIT_i$ and its unique
transform $\hat E=\Phi^i(E)$ is locally-free.
\end{defin}

In the situations in which we can apply the theorem on cohomology
and base change we can express this condition in the following
equivalent way.

\begin{rem}
If the kernel $\mathcal P$ is a locally free sheaf then $E$ is
$\Phi$-$\IT_i$ if the cohomology group $H^j (X, E \otimes \Pc_y) =
0$ vanishes for every $j\neq i$ and every $y \in Y$, where $\Pc_y
$ is the restriction of $\Pc $ to $X\times {\{y\}}$. In this case,
the fibre over $y\in Y$ of the unique transform $\hat E$ is
canonically isomorphic to $H^i(X,E\otimes \Pc_y)$.
\end{rem}

Let us consider now an elliptic curve $X$. Given a sheaf $E$ on $X$
we write its Chern character as a pair of integers $(r(E),d(E))$.
Let $\alpha$ and $\beta$ be coprime integers with $\beta>0$ and let
$Y={M}(\beta,\alpha)$ be the moduli space of stable bundles on $X$
of Chern character $(\beta,\alpha)$. The work of Atiyah \cite{A}
(see also Tu \cite{Tu}) implies that $Y$ is isomorphic to $X$; we
preserve the distinction for clarity. Let $\Pc$ be a universal
bundle on $X\times Y$, and put
$$\Phi=\Phi^{\Pc}_{X\to Y},\qquad \Psi=\Phi^{\Pc^{\dual}}_{Y\to X}.$$
As we noted above, $\Psi[1]$ is a left adjoint of $\Phi$. The
following results on the classification of Fourier-Mukai
transforms on elliptic curves (Proposition \ref{yawn} and Theorem
\ref{ellcur}) are due to Bridgeland \cite{Br}.

\begin{prop}
\label{yawn} The functor $\Phi$ is a Fourier-Mukai transform.
\end{prop}

By the Grothendieck-Riemann-Roch theorem there is a group
homomorphism $\Phi_*$ making the following diagram commute

$$\xymatrix{ D(X)  \ar[r]^(.5){\Phi} \ar[d]_{\ch} &
D(Y) \ar[d]^{\ch} \\
H^{even}(X,\mathbb{Z}) \ar[r]^(.5){\Phi_*} &
H^{even}(Y,\mathbb{Z})}$$

Since $\Phi_*$ is an invertible morphism whose inverse is
$-\Psi_*$ it follows that $\Phi_*$ is given by a matrix of the
form $A=\mat{\alpha}{\beta}{\gamma}{\delta}$ that belongs to
$SL(2,\mathbb{Z})$.

Recall that a $Y$-flat sheaf $\Pc$ on $X\times Y$ is said to be {\it
strongly simple} over $Y$ if $\Pc_y$ is simple for all $y\in Y$, and
if for any pair $y_1,y_2$ of distinct points of $Y$ and any integer
$i$ one has $\Ext^i_X(\Pc_{y_1},\Pc_{y_2})=0$. Taking into account
the isomorphism $X\simeq Y$ we have:

\begin{thm}
\label{ellcur} Let $X$ be an elliptic curve and take an element
$$A=\mat{\alpha}{\beta}{\gamma}{\delta}\in\mbox{SL}_2(\Z),$$
such that $\beta>0$. Then there exist vector bundles on $X\times
X$ that are strongly simple over both factors, and which restrict
to give bundles of Chern character $(\beta,\alpha)$ on the first
factor and $(\beta,\delta)$ on the second. For any such bundle
$\Pc$, the resulting functor $\Phi=\Phi^{\Pc}_{X\to X}$ is a
Fourier-Mukai transform, and satisfies
$$\col{r(\Phi( E))}{d(\Phi( E))}=\mat{\alpha}{\beta}{\gamma}{\delta}\col{r(E)}{d(E)},$$
for all objects $E$ of $ D(X)$.
\end{thm}

The importance of this theorem lies in the fact that it
essentially describes all the Fourier-Mukai transforms on an
elliptic curve. In order to have a complete description we must
recall the relation that exists between these transforms and
derived autoequivalences. The fundamental result, proved by Orlov
\cite{O1}, is that any derived equivalence between bounded derived
categories of smooth projective manifolds is obtained from a
Fourier-Mukai functor whose kernel is unique up to isomorphisms;
he also determined in \cite{O2} the group of derived
auto-equivalences of an abelian variety. The particular case of
elliptic curves has been considered by Hille and van den Bergh in
\cite{Hi}, these authors proved the theorem that follows.

\begin{thm}\label{auto-equivs} Let $X$ be an elliptic curve. The group of
derived auto-equivalences $\mathrm{Aut}(D(X))$ sits in the exact
sequence

$$0\to 2\,\mathbb Z \times\mathrm{Aut}(X)\ltimes\mathrm{Pic}^0(X)\to
\mathrm{Aut}(D(X))\xrightarrow{\mathrm{ch}} \mathrm{SL}(2,\mathbb
Z)\to 0$$ where $n\in\Z$ acts by means of the shift functor $[n]$,
the transform corresponding to $(f,L)\in
\mathrm{Aut}(X)\ltimes\mathrm{Pic}^0(X)$ sends $\mathcal E$ to $
f_*(L\otimes\mathcal E)$ and for any $\Phi\in\mathrm{Aut}(D(X))$ we
have $\mathrm{ch}(\Phi)=\Phi_*$.
\end{thm}

For the rest of the paper it is crucial to have a good
understanding of the action of Fourier-Mukai transforms on vector
bundles over elliptic curves. The following result, in the case of
simple sheaves, was proved by Bridgeland \cite{Br}. The extension
to indecomposable sheaves is straightforward.

\begin{prop}
\label{stable} Let $X$ be an elliptic curve and let
$\Phi=\Phi^{\Pc}_{X\to X}\colon D(X)\to D(X)$ be a Fourier-Mukai
transform. Any simple (indecomposable) sheaf $E$ on $X$ is
$\Phi$-$\WIT$ and the transform $\Hat E$ is a simple
(indecomposable) sheaf.
\end{prop}

\begin{rem}\label{serre}
We recall that a simple (indecomposable) vector bundle on an
elliptic curve $X$ is stable (semistable). In particular, a vector
bundle on $X$ is simple if and only if it is stable (see
\emph{\cite{Br}}).
\end{rem}

Now we determine the action of any Fourier-Mukai transform on
stable and semistable bundles over an elliptic curve, these are
the main results of this section.

\begin{prop}\label{presstab} Let $\Phi\colon D(X)\to D(X)$ be a
Fourier-Mukai transform on an elliptic curve $X$ with
$\Phi_*=\mat{\alpha}{\beta}{\gamma}{\delta}$, where $\beta>0$, and
let $E$ be a semistable (stable) vector bundle of Chern character
$\ch(E)=(r,d)$ over $X$.

\begin{enumerate}
\item If $\alpha\,r+\beta\,d> 0$ then $E$ is
$\Phi$-$\IT_0$ and the transform $\widehat E$ is also semistable
(stable).

\item If $\alpha\,r+\beta\,d= 0$ then $E$ is
$\Phi$-$\WIT_1$ and the transform $\widehat E$ is a torsion sheaf.

\item If $\alpha\,r+\beta\,d< 0$ then $E$ is
$\Phi$-$\IT_1$ and the transform $\widehat E$ is also semistable
(stable).
\end{enumerate} Finally, if $E$ is $\Phi$-$\WIT_i$
then $\ch(\widehat E)=(-1)^i\, \Phi_*(\ch(E))$.
\end{prop}

\begin{proof}
By Theorem \ref{ellcur} there exists a vector bundle $\Pc\to
X\times X$ that is strongly simple over both factors and such that
$\Phi=\Phi^{\Pc}_{X\to X}$. Since $E$ is semistable and
$\Pc_{|X\times\{x\}}$ is stable with
$\ch(\Pc_{|X\times\{x\}})=(\beta,\alpha)$ it follows that
$E\otimes\Pc_{|X\times\{x\}}$ is semistable and
$\deg(E\otimes\Pc_{|X\times\{x\}})=\alpha\,r+\beta\,d$. The
$\Phi$-$\IT$ statements follow now from the Riemann-Roch Theorem
and the vanishing of all sections of a semistable vector bundle of
negative degree. In these cases the preservation of stability
(semistability) follows from Proposition \ref{stable} and Remark
\ref{serre}.

Now let us assume  that $\alpha\,r+\beta\,d= 0$. Since $\alpha$ and
$\beta$ are relatively prime it follows that $r=h\,\beta$,
$d=-h\,\alpha$ where $h=\gcd(r,d)$ is the greatest common divisor of
$r$ and $d$. Therefore, the quotients $E_i$ of any Jordan-H\"{o}lder
filtration of $E$ are stable vector bundles of slope
$-\frac{\alpha}{\beta}$. It follows that $\ch(E_i)=(\beta,-\alpha)$.
On the other hand, if $\Psi[1]$ is the quasi-inverse of $\Phi$, the
skyscraper sheaf $\mathcal{O}_x$ of a point $x\in X$ is
$\Psi-\WIT_0$ with $\Psi^0(\mathcal{O}_x)=\Pc^\vee_{|X\times\{x\}}$
and $\ch(\Pc^\vee_{|X\times\{x\}})=(\beta,-\alpha)$. Therefore every
stable bundle with Chern character $(\beta,-\alpha)$ is
$\Phi-\WIT_1$ and its transform is a torsion sheaf. The second
statement follows now by considering a Jordan-H\"{o}lder filtration.
\end{proof}

In a similar way we prove:

\begin{prop}\label{presstab-adj} Let $\Psi[1]\colon D(X)\to D(X)$ be
the quasi-inverse of $\Phi$ and let $E$ be a semistable (stable)
vector bundle of Chern character $\ch(E)=(r,d)$ over $X$.

\begin{enumerate}
\item If $-\delta\,r+\beta\,d> 0$ then $E$ is
$\Psi$-$\IT_0$ and the transform $\widehat E$ is also semistable
(stable).

\item If $-\delta\,r+\beta\,d= 0$ then $E$ is
$\Psi$-$\WIT_1$ and the transform $\widehat E$ is a torsion sheaf.

\item If $-\delta\,r+\beta\,d< 0$ then $E$ is
$\Psi$-$\IT_1$ and the transform $\widehat E$ is also semistable
(stable).
\end{enumerate} Finally, if $E$ is $\Psi$-$\WIT_i$
then $\ch(\widehat E)=(-1)^{i+1}\, (\Phi_*)^{-1}(\ch(E))$.
\end{prop}

\section{Fourier-Mukai transforms for coherent systems}
\label{trans-coh-sys}

\subsection{Coherent systems}

A coherent system of type $(r,d,k)$ on a smooth projective curve $X$
is by definition a pair $(E,V)$ consisting of a vector bundle $E$ of
rank $r$ and degree $d$ over $X$ and a vector subspace $V \subset
H^0(X,E)$ of dimension $k$. A morphism $f\colon
(E^\prime,V^\prime)\to (E,V)$ of coherent systems is a homomorphism
of vector bundles $f\colon E^\prime\to E$ such that
$f(V^\prime)\subset V$. If $E^\prime$ is a subbundle of $E$ then we
say that $(E^\prime,V^\prime)$ is a coherent subsystem of $(E,V)$.
With this definition, coherent systems on $X$ form an additive
category $S(X)$ (see \cite{LP} \S 4.1).

For any real number $\alpha$, the  $\alpha$-slope of a coherent
system $(E,V)$ of type $(r,d,k)$ is defined by
$$
\mu_{\alpha}(E,V) = \frac{d}{r} + \alpha \frac{k}{r}.
$$
A coherent system $(E,V)$ is called $\alpha$-stable ({
$\alpha$-semistable}) if
$$
\mu_{\alpha}(E',V') < \mu_{\alpha}(E,V) \ \ (\mu_{\alpha}(E',V')
\le \mu_{\alpha}(E,V))
$$
for every proper coherent subsystem $(E',V')$ of $(E,V)$.

Following \cite{KN}, by replacing a coherent system $(E,V)$ by its
evaluation map $V\otimes\mathcal{O}_X \to E$ we may regard
coherent systems as forming a subcategory of the abelian category
$\mathcal{C}(X)$ whose objects are arbitrary sheaf maps
$\varphi\colon V\otimes\mathcal{O}_X \to \mathcal{E}$ where $V$ is
a finite dimensional vector space and $\mathcal{E}$ is any
coherent sheaf. A morphism in $\mathcal{C}(X)$ from
$\varphi_1\colon V_1\otimes\mathcal{O}_X \to \mathcal{E}_1$ to
$\varphi_2\colon V_2\otimes\mathcal{O}_X \to \mathcal{E}_2$ is
given by a linear map $f\colon V_1\to V_2$ and a sheaf map
$g:\mathcal{E}_1\to \mathcal{E}_2$ such that the following diagram
commutes

$$\xymatrix{ V_1\otimes\mathcal{O}_X  \ar[r]^(.6){\varphi_1} \ar[d]_{f\otimes 1} &
{\mathcal{E}_1} \ar[d]^{g} \\
V_2\otimes\mathcal{O}_X  \ar[r]^(.6){\varphi_2} & {\mathcal{E}_2}.
}$$ With this definition, the category $S(X)$ of coherent systems
is a full subcategory of $\mathcal C(X)$.

\begin{rem}\label{repr}
An object $\varphi\colon V\otimes\mathcal{O}_X \to \mathcal{E}$ of
$\mathcal{C}(X)$ represents a coherent system if and only if
$\mathcal{E}$ is a vector bundle and the induced map
$H^0(\varphi)\colon V\to H^0(X,\mathcal{E})$ is injective.
Considering the natural exact sequence $$0\to N\xrightarrow{j}
V\otimes\mathcal{O}_X \xrightarrow{\varphi} \mathcal E\to F\to 0$$
where $N=\Ker\varphi$, $F=\Coker\varphi$, the injectivity of
$H^0(\varphi)$ is equivalent to $H^0(X,N)=0$.
\end{rem}

The notion of $\alpha$-(semi)stability of coherent systems can be
extended to $\mathcal{C}(X)$. This extension does not introduce
new semistable objects. That is, an object $\varphi$ of
$\mathcal{C}(X)$ is $\alpha$-(semi)stable if and only if it is an
$\alpha$-(semi)stable coherent system (see \cite{KN} \S 2).

The full subcategory $S_{\alpha,\mu}(X)$ of $\mathcal{C}(X)$
consisting of $\alpha$-semistable coherent systems with fixed
$\alpha$-slope $\mu$ is a Noetherian and Artinian abelian category
whose simple objects are precisely the $\alpha$-stable coherent
systems (see \cite{KN,RV}).

The $\alpha$-stable coherent systems of type $(r,d,k)$ on $X$ form
a quasiprojective moduli space that we denote by
$G(\alpha;r,d,k)$. In the particular case where $X$ is an elliptic
curve, Lange and Newstead \cite{LN} have obtained a complete
description of these moduli spaces. In order to state their
results, we first recall the definition of the Brill-Noether
number $\beta(r,d,k)$ \cite[Definition 2.7]{BGMN}. For an elliptic
curve, this is independent of $r$ and we write
$$
\beta(d,k):=\beta(r,d,k)=k(d-k)+1.
$$
As before, $M(r,d)$ denotes the moduli space of stable bundles of
rank $r$ and degree $d$ on $X$. Note that, if $k\ge1$ and
$\alpha\le0$, there do not exist $\alpha$-stable coherent systems
of type $(r,d,k)$ (see \cite[Section 2.1]{BGMN}). The main results
of \cite{LN} can now be summarized as follows.

\noindent{\bf Theorem.}\begin{em} Let $X$ be an elliptic curve and
suppose $r\ge1$, $k\ge0$. Then
\begin{itemize}
\item[(i)] if $G(\alpha;r,d,k)\ne\emptyset$, it is smooth and irreducible
of dimension $\beta(d,k)$;
\item[(ii)] $G(\alpha;r,d,0)\simeq M(r,d)$ for all $\alpha$; in particular
it is non-empty if and only if $\gcd(r,d)=1$;
\item[(iii)] for $\alpha>0$ and $k\ge1$, $G(\alpha;1,d,k)$ is independent of
$\alpha$ and is non-empty if and only if either $d=0$, $k=1$ or
$k\le d$;
\item[(iv)] for $\alpha>0$, $r\ge2$ and $k\ge1$, $G(\alpha;r,d,k)\ne\emptyset$
if and only if $(r-k)\alpha<d$ and either $k<d$ or $k=d$ and
$\gcd(r,d)=1$.
\end{itemize}\end{em}

\subsection{Fourier-Mukai transforms}\label{fm-coh-sys}

As a first step to transform coherent systems on an elliptic curve
$X$ be begin by studying the Fourier-Mukai transforms $\Phi\colon
D(X)\to D(X)$ that act on the abelian category $\mathcal{C}(X)$. A
natural requirement for such a transform $\Phi$ is that
$\mathcal{O}_X$ should be $\Phi$-$\WIT_i$ and
$\Phi^i(\mathcal{O}_X)=\mathcal{O}_X$.

\begin{prop}

Let $a$ be a positive integer. There exists a Fourier-Mukai
transform $\Phi_a\colon D(X)\to D(X)$, unique up to composition
with automorphisms of $X$, such that $\mathcal{O}_X$ is
$\Phi_a$-$\IT_0$, $\Phi_a^0(\mathcal{O}_X)=\mathcal{O}_X$ and
$(\Phi_a)_*=\mat{1}{a}{0}{1}$.

If $\Psi_a[1]\colon D(X)\to D(X)$ is the quasi-inverse of $\Phi_a$
then
 $\mathcal{O}_X$ is $\Psi_a$-$\IT_1$,
$\Psi_a^1(\mathcal{O}_X)=\mathcal{O}_X$ and
$(\Psi_a)_*=\mat{-1}{a}{0}{-1}$.
\end{prop}

\begin{proof}
Let $\Phi\colon D(X)\to D(X)$ be a Fourier-Mukai transform such that
$\Phi_*(\ch(\mathcal{O}_X))=\ch(\mathcal{O}_X)$. Then
$\Phi_*=\mat{\alpha}{\beta}{\gamma}{\delta}$ has
$\ch(\mathcal{O}_X)=(1,0)$ as an eigenvector, therefore $\alpha=1$
and $\gamma=0$. Moreover, since $\Phi_*\in SL(2,\mathbb{Z})$ we must
have $\delta=1$.  If we set $\beta=a$, we get
${\Phi}_*=\mat{1}{a}{0}{1}$. By Theorem \ref{ellcur} there exist
vector bundles $\Pc$ on $X\times X$, strongly simple over both
factors, such that $\Phi_*=(\Phi^\Pc_{X\to X})_*$. Moreover, by
Theorem \ref{auto-equivs}, $\Phi$ is determined up to composition
with the transforms $f_*(L\otimes (-))$ where $(f,L)\in
\mathrm{Aut}(X)\ltimes\mathrm{Pic}^0(X)$. Therefore, we can write
$\Phi=f_*(L\otimes (-))\circ\Phi^{\Pc}_{X\to X}$ and
$\Phi(\mathcal{O}_X)=f_*(L\otimes (\Phi^{\Pc}_{X\to
X}(\mathcal{O}_X))) $. The condition
$\Phi^0(\mathcal{O}_X)=\mathcal{O}_X$ uniquely determines $L$ and
leaves $f\in \mathrm{Aut}(X)$ undetermined.

\end{proof}

\begin{defin}
The Fourier-Mukai transforms $\Phi_a$, $\Psi_a$ define functors
$\Phi_a^0\colon \mathcal{C}(X)\to \mathcal{C}(X)$, $\Psi_a^1\colon
\mathcal{C}(X)\to \mathcal{C}(X)$ that send the object
$\varphi\colon V\otimes\mathcal{O}_X \to \mathcal{E}$ to
\begin{align*}&\Phi_a^0(\varphi)\colon
V\otimes\mathcal{O}_X \to \Phi_a^0(\mathcal{E})\\
\\
&\Psi_a^1(\varphi)\colon V\otimes\mathcal{O}_X \to
\Psi_a^1(\mathcal{E}),
\end{align*}respectively.
\end{defin}

Let us recall that for any integral functor $\Phi\colon D(X)\to
D(Y)$ and any sheaf $\mathcal{E}$ there is a spectral sequence
$$E_2^{p,q}=H^p(X,\Phi^q(\mathcal{E}))\Rightarrow
H^{p+q}(Y,\Phi(\mathcal{E}))$$ that is obtained from the spectral
sequence for the composition of two derived functors.

\begin{remark}\label{part-case}
In the particular case of an integral functor $\Phi\colon D(X)\to
D(X)$ on an elliptic curve $X$, the spectral sequence collapses at
the $2^\mathrm{nd}$ term and for every integer $n$ we get an exact
sequence $$0\to H^1(X,\Phi^{n-1}(\mathcal{E}))\to
H^n(X,\Phi(\mathcal{E}))\to H^0(X,\Phi^n(\mathcal{E}))\to 0$$
\end{remark}

\begin{lemma}\label{suc-esp} Let $X$ be an elliptic curve and let us consider
the Fourier-Mukai transform $\Phi_a\colon D(X)\to D(X)$ and its
quasi-inverse $\Psi_a[1] \colon D(X)\to D(X)$.

\begin{enumerate}
\item There is an exact sequence $$0\to H^1(X,\Phi^{0}_a(\mathcal{E}))\to
H^1(X,\mathcal{E})\to H^0(X,\Phi^1_a(\mathcal{E}))\to 0$$ together
with the information $H^0(X,\Phi^{0}_a(\mathcal{E}))\simeq
H^0(X,E)$, $H^1(X,\Phi^{1}_a(\mathcal{E}))=0$.

\item There is an exact sequence $$0\to H^1(X,\Psi^{0}_a(\mathcal{E}))\to
H^0(X,\mathcal{E})\to H^0(X,\Psi^1_a(\mathcal{E}))\to 0$$ together
with the information $H^0(X,\Psi^{0}_a(\mathcal{E}))=0$,
$H^1(X,\Psi^{1}_a(\mathcal{E}))\simeq H^1(X,E)$.
\end{enumerate}
\end{lemma}

\begin{proof} Since $\Phi_a$ is an equivalence and
$\mathcal{O}_X$ is $\Phi_a-\WIT_0$ we have
$H^n(X,\Phi_a(\mathcal{E}))\simeq H^n(X,\mathcal{E})$. In a
similar way we get $H^n(X,\Psi_a(\mathcal{E}))\simeq
H^{n-1}(X,\mathcal{E})$. According to Remark \ref{part-case} this
is enough.
\end{proof}

\begin{prop}\label{fm-cs-I}
For any coherent system $\varphi\colon V\otimes\mathcal{O}_X \to
E$ its $\Phi_a^0$-transform $$\Phi_a^0(\varphi)\colon
V\otimes\mathcal{O}_X \to \Phi_a^0(E)$$ is a coherent system.
Therefore, the functor $\Phi_a^0\colon\mathcal{C}(X)\to
\mathcal{C}(X)$ leaves the subcategory $S(X)$ of coherent systems
stable and we get an induced functor
$$\Phi_a^0\colon S(X)\to S(X).$$
\end{prop}

\begin{proof} Applying the transform $\Phi_a^0$ to the natural exact sequence
$$0\to N\xrightarrow{j}
V\otimes\mathcal{O}_X \xrightarrow{\varphi}  E\to F\to 0$$ we get
the exact sequence $0\to\Phi_a^0(N)\xrightarrow{\Phi_a^0(j)}
V\otimes\mathcal{O}_X \xrightarrow{\Phi_a^0(\varphi)}
\Phi_a^0(E)$. Therefore, according to Remark \ref{part-case}, the
transform is a coherent system if and only if
$H^0(X,\Phi_a^0(N))=0$ and $\Phi_a^0(E)$ is locally free.

By Lemma \ref{suc-esp} we have $H^0(X,\Phi_a^0(N)=H^0(X,N)$ and
Remark \ref{part-case} implies that $H^0(X,N)=0$ since
$\varphi\colon V\otimes\mathcal{O}_X \to E$ is a coherent system. On
the other hand, since we are on a curve and $\Phi_a$ is a
Fourier-Mukai transform it follows that $\Phi_a^0(E)$ is
$\Psi_a-\WIT_1$ and its transform is $\Phi_a-IT_0$, hence
$\Phi_a^0(E)$ is locally free.
\end{proof}

\begin{prop}\label{psi-it1}
If $\varphi\colon V\otimes\mathcal{O}_X \to E$ is a coherent
system such that $E$ is $\Psi_a-\IT_1$ then its
$\Psi_a^1$-transform
$$\Psi_a^1(\varphi)\colon V\otimes\mathcal{O}_X \to \Psi_a^1(E)$$
is a coherent system.
\end{prop}

\begin{proof}
The natural exact sequence $$0\to N\xrightarrow{j}
V\otimes\mathcal{O}_X \xrightarrow{\varphi}  E\xrightarrow{\pi}
F\to 0$$ splits into the short exact sequences
\begin{equation}\label{A} 0\to N \xrightarrow{j} V\otimes\mathcal{O}_X
\xrightarrow{\bar\varphi}\Ima\varphi\to 0\end{equation}
\begin{equation}\label{B}0\to \Ima\varphi\xrightarrow{i} E\xrightarrow{\pi}
F\to 0 \end{equation} The transform $\Phi_a^1(\varphi)\colon
V\otimes\mathcal{O}_X \to \Phi_a^1(E)$ is a coherent system if and
only if $H^0(\Psi_a^1(\varphi))$ is injective. By (\ref{A}) and
(\ref{B}) $\varphi=i\circ\bar\varphi$, thus
$H^0(\Psi_a^1(\varphi))=H^0(\Psi_a^1(i))\circ
H^0(\Psi_a^1(\bar\varphi))$.

If we apply the Fourier-Mukai transform $\Psi_a$ to (\ref{A}) and
(\ref{B}) and we take into account that $\mathcal{O}_X$ and $E$
are $\Psi_a-\IT_1$ we get the exact sequences

\begin{equation}\label{I} 0\to
\Psi_a^1(N)\xrightarrow{\Psi_a^1(j)} V\otimes\mathcal{O}_X
\xrightarrow{\Psi_a^1(\bar\varphi)}\Psi^1_a(\Ima\varphi)\to
0\end{equation}

\begin{equation}\label{II}0\to \Psi_a^0(F)\to
\Psi_a^1(\Ima\varphi)\xrightarrow{\Psi_a^1(i)}
\Psi_a^1(E)\xrightarrow{\Psi_a^1(\pi)} \Psi_a^1(F)\to 0
\end{equation} Splitting (\ref{II}) we get the short exact
sequences

\begin{equation}\label{IIA}
0\to \Psi_a^0(F) \to
\Psi_a^1(\Ima\varphi)\xrightarrow{f}\Ima\Psi_a^1(i)\to 0
\end{equation}

\begin{equation}\label{IIB}
0\to \Ima\Psi_a^1(i)\xrightarrow{g} \Psi_a^1(E)\to \Psi_a^1(F)
\to 0
\end{equation} where $\Psi_a^1(i)=g\circ f$.

Now if we take cohomology in (\ref{I}), (\ref{IIA}), (\ref{IIB}) and
bear in mind Lemma \ref{suc-esp}, we get that
$H^0(\Psi_a^1(\bar\varphi))$ and $H^0(\Psi_a^1(i))$ are injective.
This finishes the proof.
\end{proof}

\begin{remark}
The functor $\Psi_a^1\colon\mathcal{C}(X)\to \mathcal{C}(X)$
induces a functor $$\Psi_a^1\colon S^1_a(X)\to S(X)$$ where
$S^1_a(X)$ is the full subcategory of $S(X)$ whose objects are
those coherent systems $\varphi\colon V\otimes\mathcal{O}_X \to E$
such that $E$ is $\Psi_a-\IT_1$.
\end{remark}

\section{Preservation of stability}\label{pres-stab}

Let us recall that $\alpha$-stable coherent systems do exist only
for $\alpha>0$ (see \cite{BGMN}). Moreover, the $\alpha$-range is
divided into open intervals determined by a finite number of
critical values
$$0=\alpha_0<\alpha_1<\cdots<\alpha_L$$ such that the
moduli spaces for any two values of $\alpha$ in the interval
$(\alpha_i,\alpha_{i+1})$ coincide; if $k\geq r$ this is also true
for the interval $(\alpha_L,\infty)$ (see \cite[Propositions 4.2 \&
4.6]{BGMN}).

In the rest of this section we assume that $d\neq 0$, $k>0$. Then,
according to \cite[Proposition 3.2 \& Lemma 4.1]{LN}, we have

\begin{prop}\label{comp} Let $(E,V)$ be an $\alpha$-stable coherent system
of type $(r,d,k)$ on an elliptic curve $X$. Then every
indecomposable direct summand of $E$ is of positive degree.
\end{prop}

\subsection{Small $\alpha$}\label{sm-alpha} For a projective curve $X$ we denote by $G_0(r,d,k)$ the moduli
space of $\alpha$-stable coherent systems of type $(r,d,k)$ with
$0<\alpha<\alpha_1$ and $\alpha_1$ is the first critical value.

The following is well known (see \cite[Remark 1.4 (ii)]{RV}).

\begin{prop}\label{esta} A coherent system $(E,V)$ of type $(r,d,k)$ is
$\alpha$-stable, with $0<\alpha<\alpha_1$, if and only if $E$ is
semistable and one has $$\frac{k^\prime}{r^\prime}< \frac{k}{r},$$
for all coherent subsystems $(E^\prime,V^\prime)$ of type
$(r^\prime,d^\prime,k^\prime)$ with $0\neq E^\prime\neq E$ and
$\mu(E^\prime)=\mu(E)$.
\end{prop}

In the rest of this subsection $\alpha$-stability refers to a
positive $\alpha$ which is less than the first critical value.

Let us assume now that $X$ is an elliptic curve.

\begin{thm}\label{small-fi}
The Fourier-Mukai transform $\Phi_a$ induces a map
$$\Phi_a^0\colon G_0(r,d,k)\to G_0(r+ad,d,k)$$
\end{thm}

\begin{proof}
Given a coherent system $\varphi\colon V\otimes\mathcal{O}_X \to
E\in G_0(r,d,k)$ by Propositions \ref{esta} and \ref{comp} $E$ is
semistable of positive degree. Hence, $r+ad>0$ and Proposition
\ref{presstab} implies that $E$ is $\Phi_a-\IT_0$ and
$\Phi_a^0(E)$ is semistable with Chern character
$\ch(\Phi_a^0(E))=(r+ad,d)$. Therefore, we can use Proposition
\ref{esta} in order to prove that $\Phi^0_a(\varphi)\colon
V\otimes\mathcal{O}_X \to \Phi^0_a(E)\in G_0(r,d,k)$. If it were
not the case, then we would have a coherent subsystem $(\widehat
E^\prime,V^\prime)$ of type $(\hat r^\prime, \hat
d^\prime,k^\prime)$ with $\mu(\widehat E^\prime)=\mu(\Phi_a^0(E))$
and
\begin{equation}\label{d1}
\frac{k}{r+ad}-\frac{k^\prime}{\hat r^\prime}\leq 0.
\end{equation}

The system and subsystem  fit into the commutative diagram

$$\xymatrix{ 0\ar[r] & {\widehat E^\prime}\ar[r] & {\Phi_a^0(E)}\ar[r] & {\widehat
E^{\prime\prime}}\ar[r] & 0\\
0\ar[r] & V^\prime\otimes\mathcal{O}_X  \ar[r] \ar[u]  &
V\otimes\mathcal{O}_X  \ar[r]\ar[u] & V/V^\prime\otimes\mathcal{O}_X
\ar[r]\ar[u] & 0 }$$

Since $\Phi_a^0(E)$ is semistable, it follows that $\widehat
E^\prime$ and $\widehat E^{\prime\prime}$ are semistable with
$\mu(\widehat E^{\prime})=\mu(\widehat E^{\prime\prime})$. Thus
$\widehat E^{\prime}$, $\widehat E^{\prime\prime}$ are
$\Psi_a-\IT_1$ and by applying $\Psi_a$ we get a subsystem
$$\xymatrix{ 0\ar[r] & {\Psi_a^1(\widehat E^\prime) }\ar[r] & {E}\\
0\ar[r] & V^\prime\otimes\mathcal{O}_X  \ar[r] \ar[u] &
V\otimes\mathcal{O}_X  \ar[u] }$$with $\ch(\Psi_a^1(\widehat
E^\prime))=(\hat r^\prime-a\hat d,\hat d)$. The stability of
$V\otimes\mathcal{O}_X \to {E}$ gives
\begin{equation}\label{d2}\frac{k}{r}-\frac{k^\prime}{\hat r^\prime-a\hat
d}>0
\end{equation}

We can write the inequalities (\ref{d1}), (\ref{d2}) as

\begin{align}
\label{d3}&\frac{k}{r(1+a\mu(E))}-\frac{k^\prime}{\hat r^\prime}\leq 0\\
\label{d4}&\frac{k}{r}-\frac{k^\prime}{\hat
r^\prime(1-a\mu(\widehat E^\prime))}> 0
\end{align} On the other hand, we have
$\mu(\Phi_a^0(E))=\frac{\mu(E)}{1+a\mu(E)}$. Thus,
$1-a\mu(\widehat E^\prime)=\frac{1}{1+a\mu(E)}$ and the inequality
(\ref{d4}) is equivalent to
$$\frac{k}{r(1+a\mu(E))}-\frac{k^\prime}{\hat r^\prime}> 0$$ which
contradicts (\ref{d3}).
\end{proof}

In a similar way we prove:

\begin{thm}\label{small-psi}
The Fourier-Mukai transform $\Psi_a$ induces a map
$$\Psi_a^1\colon G_0(r+ad,d,k)\to G_0(r,d,k).$$
\end{thm}

As a consequence of Theorems \ref{small-fi}, \ref{small-psi}, and
the fact that $\Psi_a[1]$ and $\Phi_a$ are quasi-inverses we get
the following result.

\begin{thm}\label{sm-iso}
The Fourier-Mukai transform $\Phi_a$ induces an isomorphism of
moduli spaces
$$\Phi_a^0\colon G_0(r,d,k)\to G_0(r+ad,d,k)$$ Therefore, the
isomorphism type of $G_0(r,d,k)$  depends only on the class
$[r]\in\mathbb Z/d\,\mathbb Z$.
\end{thm}

\subsection{Large $\alpha$}\label{lg-alpha}

In this subsection we suppose that $0<k<r$. Under this assumption,
Lange and Newstead proved in \cite[Theorem 5.2]{LN} that for an
elliptic curve the moduli space $G(\alpha;r,d,k)$ is non empty if
and only if $0<\alpha<\frac{d}{r-k}$ and either $k<d$ or $k=d$ and
$\gcd(r,d)=1$. Moreover, in this case the largest critical value
$\alpha_L$ verifies $\alpha_L<\frac{d}{r-k}$.

For a projective curve $X$ we denote by $G_L(r,d,k)$ the moduli
space of $\alpha$-stable coherent systems of type $(r,d,k)$ with
$\alpha_L<\alpha<\frac{d}{r-k}$. The description of this moduli
space by means of BGN extensions (see {\cite{BGN}}) has been carried
out in \cite{BG} where we refer for details (see also \cite[Remark
5.5]{BGMN}). Now we recall the main concepts used in these
references.

\begin{defin}\label{bgn}
A BGN extension of type $(r,d,k)$ is an extension of vector
bundles
$$0\to \mathcal{O}_X^k\to E\to F\to 0$$ which satisfies the
following conditions
\begin{enumerate}
\item $\rk E=r>k$
\item $\deg E=d>0$
\item $H^0(X,F^\vee)=0$
\item If $(e_1,\ldots,e_k)\in
\mathrm{Ext}_X^1(F,\mathcal{O}_X^k)\simeq H^1(X,F^\vee)^{\oplus
k}$ denotes the class of the extension, then $e_1,\ldots,e_k$ are
linearly independent as vectors in $H^1(X,F^\vee)$.
\end{enumerate}
\end{defin}

The following result corresponds to \cite[Proposition 4.1]{BG}.

\begin{prop}\label{est-bgn} Let $(E,V)$ be an $\alpha$-semistable
coherent system of type $(r,d,k)$ with
$\alpha_L<\alpha<\frac{d}{r-k}$. Then the evaluation map of
$(E,V)$ defines a BGN extension $$0\to \mathcal{O}_X^k\to E\to
F\to 0$$ with $F$ semistable.
\end{prop}

On the converse direction we prove now:

\begin{prop}\label{esta-bgn} On a
projective curve $X$ a BGN extension of type $(r,d,k)$
$$0\to \mathcal{O}_X^k\to E\to F\to 0$$ defines an
$\alpha$-stable coherent system, with
$\alpha_L<\alpha<\frac{d}{r-k}$, if and only if $F$ is semistable
and one has $$\frac{k^\prime}{r^\prime}> \frac{k}{r},$$ for all
subextensions $$0\to \mathcal{O}_X^{k^\prime}\to E^\prime\to
F^\prime\to 0$$ of type $(r^\prime,d^\prime,k^\prime)$ with
 $\mu(F^\prime)=\mu(F)$.
\end{prop}

\begin{proof}
Since $\alpha_L<\alpha<\frac{d}{r-k}=\mu(F)$ we can write
$\alpha=\mu(F)-\epsilon$, with $\epsilon>0$ and small enough. The
$\alpha$-slope $\mu_\alpha(E,k)$ of $\mathcal{O}_X^k\to E$ can be
expressed as $\mu_\alpha(E,k)=\mu(F)-\epsilon\,\frac{k}{r}$.
Given a coherent subsystem $\mathcal{O}_X^{k^\prime}\to E^\prime$
we have a commutative diagram
$$\xymatrix{ & 0 & 0 & 0 & \\
& {\mathcal{O}_X^{k-k^\prime}} \ar[u]\ar[r] & E^{\prime\prime}
\ar[u]\ar[r] & F^{\prime\prime} \ar[u] &\\
0\ar[r] & {\mathcal{O}_X^k}\ar[r]\ar[u] & E \ar[r]\ar[u] &
F\ar[r]\ar[u] & 0\\
0\ar[r] & {\mathcal{O}_X^{k^\prime}}\ar[r]\ar[u] & E^\prime
\ar[r]\ar[u] &
F^\prime\ar[r]\ar[u] & 0\\
& 0\ar[u] & 0\ar[u] & K\ar[u] &\\
& & & 0\ar[u]& }$$whose rows and columns are exact, and we have
the additional information given by the associated $\Ker-\Coker$
exact sequence $0\to K\to {\mathcal{O}_X^{k-k^\prime}}\to
E^{\prime\prime}\to F^{\prime\prime}$.

We have
$\mu_\alpha(E^\prime,k^\prime)=\mu(E^\prime)+(\mu(F)-\epsilon)\,\frac{k}{r}$,
and the difference function is
$$\Delta_\alpha(E^\prime,k^\prime)=\mu_\alpha(E,k)-\mu_\alpha(E^\prime,k^\prime)
= \left(\frac{r^\prime-k^\prime}{r^\prime}\right)
\mu(F)-\mu(E^\prime)+
\epsilon\,\left(\frac{k^\prime}{r^\prime}-\frac{k}{r}\right).$$

Note that $\rk F^\prime=0$ if and only if $F^\prime=0$. In fact, if
$F^\prime$ were a torsion sheaf then $K=0$ since it is a subsheaf of
the locally free sheaf ${\mathcal{O}_X^{k-k^\prime}}$. Thus
$F^\prime$ would be a subsheaf of the locally free sheaf $F$, which
is impossible. On the other hand, if $F^\prime=0$ we have
$\Delta_\alpha(E^\prime,k^\prime)=\epsilon\,(\frac{r-k}{r})>0$.
Hence, this type of coherent subsystems are not destabilizing, and
we can assume from now on, without loss of generality, that $\rk
F^\prime\neq 0$.

Therefore we can write
$\mu(E^\prime)=(\frac{r^\prime-k^\prime}{r^\prime})\,\mu(F^\prime)$
and the difference function is
\begin{equation}\label{diff}
\Delta_\alpha(E^\prime,k^\prime)=\left(\frac{r^\prime-k^\prime}{r^\prime}
\right) (\mu(F)-\mu(F^\prime))+
\epsilon\,\left(\frac{k^\prime}{r^\prime}-\frac{k}{r}\right).
\end{equation}

Let us analyze the behavior of the term
$\delta(F^\prime)=\mu(F)-\mu(F^\prime)$. In order to do this let us
split the exact sequence $0\to K\to F^\prime\to F\to
F^{\prime\prime}\to 0$ into the short exact sequences
\begin{align*}
& 0\to K\to F^\prime\to I\to 0\\
& 0\to I\to F\to F^{\prime\prime}\to 0
\end{align*}

Since $K$ is locally free, its rank is zero if and only if $K=0$.
In that case $F^\prime$ is a subsheaf of $F$ and
$\delta(F^\prime)\geq 0$ due to the semistability of $F$.

In the same way $\rk I=0$ if and only if $I=0$. If we assume
$I=0$, then $F^\prime\simeq K$ and
$\mu(F^\prime)=\mu(K)\leq\mu(\mathcal{O}_X^{k-k^\prime})=0$ since
$K$ is a subsheaf of $\mathcal{O}_X^{k-k^\prime}$. On the other
hand $\mu(F)>0$ and we get $\delta(F^\prime)>0$.

In case that $r_K=\rk K\neq 0$, $r_I=\rk I\neq 0$ we have
$\mu(F)=\frac{r_K}{r_K+r_I}\,\mu(K)+\frac{r_I}{r_K+r_I}\,\mu(I)$,
thus $$\delta(F^\prime)=\frac{r_K}{r_K+r_I}\,\mu(F)+
+\frac{r_I}{r_K+r_I}\,(\mu(F)-\mu(I))-\frac{r_K}{r_K+r_I}\mu(K)>0$$
since $\mu(F)>0$, $\mu(F)-\mu(I)\geq 0$ and $\mu(K)\geq 0$.

Due to the boundedness of the rational numbers that appear in the
difference formula (\ref{diff}), and since we can chose $\epsilon$
as small as we wish, the only coherent subsystems that can be
destabilizing are those with $\delta(F^\prime)=0$. According to
the analysis that we have made above, this is possible if and only
if $K=0$, and this is precisely the case corresponding to a
subextension.
\end{proof}

\begin{rem}\label{st-quot}
Any BGN extension $0\to \mathcal{O}_X^k\to E\to F\to 0$ in which the
quotient $F$ is stable gives rise to an $\alpha$-stable coherent
system, with $\alpha_L<\alpha<\frac{d}{r-k}$ (see \cite[Proposition
4.2]{BG}).
\end{rem}

Let us denote by $\mathcal{BGN}(r,d,k)$, $\mathcal{BGN}^s(r,d,k)$
the families of BGN extension classes of type $(r,d,k)$ on $X$ in
which the quotient is semistable or stable, respectively. Notice
that Remark \ref{st-quot} and Proposition \ref{est-bgn} establish
the following inclusions $$\mathcal{BGN}^s(r,d,k)\hookrightarrow
G_L(r,d,k)\hookrightarrow \mathcal{BGN}(r,d,k).$$

Let us assume now that $X$ is an elliptic curve.

\begin{prop}\label{bgn-phi}
The Fourier-Mukai transform $\Phi_a$ induces a map
$$\Phi_a^0\colon \mathcal{BGN}(r,d,k)\to \mathcal{BGN}(r+ad,d,k)$$
by sending a BGN extension $0\to \mathcal{O}_X^k\to E\to F\to 0$ to
$$0\to \mathcal{O}_X^k\to \Phi_a^0(E)\to \Phi_a^0(F)\to 0.$$

Moreover, $\Phi_a^0$ restricts to a map $$\Phi_a^0\colon
\mathcal{BGN}^s(r,d,k)\to \mathcal{BGN}^s(r+ad,d,k).$$
\end{prop}

\begin{proof}
Due to Proposition \ref{presstab} $F$ is  $\Phi_a-\IT_0$ and
$\Phi_a^0(F)$ is semistable. Since $\mathcal{O}_X$ is also
$\Phi_a-\IT_0$, it follows that we have an exact sequence $0\to
\mathcal{O}_X^k\to \Phi_a^0(E)\to \Phi_a^0(F)\to 0$. By Lemma
\ref{suc-esp} we have $H^0(X,\Phi_a^0(F)^\vee)=0$. On the other
hand, since $\Phi_a$ is an equivalence we get an isomorphism
$$\phi_a\colon \mathrm{Ext}_X^1 (F,\mathcal{O}_X)^{\oplus
k}\simeq\mathrm{Ext}_X^1(F,\mathcal{O}_X^k)=\mathrm{Ext}_X^1
(\Phi_a^0(F),\mathcal{O}_X^k)\simeq \mathrm{Ext}_X^1
(\Phi_a^0(F),\mathcal{O}_X)^{\oplus k}$$ Therefore if
$(e_1,\ldots,e_k)$ represents the class of the original extension,
then the transformed extension is represented by
$(\phi_a(e_1),\ldots,\phi_a(e_k))$ and these vectors continue to be
linearly independent.

Finally, if $F$ is stable then it follows from Proposition
\ref{presstab} that $\Phi_a^0(F)$ is stable, this proves the last
statement.
\end{proof}

In a similar way we prove:

\begin{prop}\label{bgn-psi}
The Fourier-Mukai transform $\Psi_a$ induces a map
$$\Psi_a^1\colon \mathcal{BGN}(r+ad,d,k)\to \mathcal{BGN}(r,d,k)$$
by sending a BGN extension $0\to \mathcal{O}_X^k\to E\to F\to 0$ to
$$0\to \mathcal{O}_X^k\to \Psi_a^1(E)\to \Psi_a^1(F)\to 0.$$

Moreover, $\Psi_a^1$ restricts to a map $$\Psi_a^1\colon
\mathcal{BGN}^s(r+ad,d,k)\to \mathcal{BGN}^s(r,d,k).$$
\end{prop}

Since $\Phi_a$ and $\Psi_a[1]$ are quasi-inverses, by Propositions
\ref{bgn-phi}, \ref{bgn-psi} we have:

\begin{thm}\label{iso-bgn}
The Fourier-Mukai transform $\Phi_a$ induces an isomorphism
$$\Phi_a^0\colon \mathcal{BGN}(r,d,k)\to \mathcal{BGN}(r+ad,d,k)$$
that restricts to an isomorphism
$$\Phi_a^0\colon \mathcal{BGN}^s(r,d,k)\to
\mathcal{BGN}^s(r+ad,d,k).$$
\end{thm}

\begin{thm}\label{large-fi}
The Fourier-Mukai transform $\Phi_a$ induces a map
$$\Phi_a^0\colon G_L(r,d,k)\to G_L(r+ad,d,k)$$
\end{thm}

\begin{proof} Given a coherent system in $G_L(r,d,k)$ we know by
Proposition \ref{est-bgn} that it defines an extension $0\to
\mathcal{O}_X^k\to E\to F\to 0$ that belongs to
$\mathcal{BGN}(r,d,k)$. By Proposition \ref{bgn-phi} the
transformed extension $0\to \mathcal{O}_X^k\to \Phi_a^0(E)\to
\Phi_a^0(F)\to 0$ belongs to $\mathcal{BGN}(r+ad,d,k)$. Therefore,
in order to prove that its is $\alpha$-stable we can use
Proposition \ref{esta-bgn}.

If it were not $\alpha$-stable then there would exist a
subextension of type $(\hat r^\prime,\hat d^\prime,k^\prime)$

$$\xymatrix{ 0\ar[r] & {\mathcal{O}_X^k}\ar[r] & \Phi_a^0(E)
\ar[r] &
\Phi_a^0(F)\ar[r] & 0\\
0\ar[r] & {\mathcal{O}_X^{k^\prime}}\ar[r]\ar[u] & {\widehat
E}^\prime \ar[r]\ar[u] &
{\widehat F}^\prime\ar[r]\ar[u] & 0\\
& 0\ar[u] & 0\ar[u] & 0\ar[u] & }$$ with $\mu(\widehat
F^\prime)=\mu(\Phi_a^0(F))$ and such that

\begin{equation}\label{bgn1}
\frac{k^\prime}{\hat r^\prime}-\frac{k}{r+ad}\leq 0
\end{equation}

Using arguments similar to those in Theorem \ref{small-fi} one
shows that after applying $\Psi_a$ we obtain a subextension of
type $(\hat r^\prime-a\hat d^\prime,\hat d^\prime,k^\prime)$ of
the original BGN extension

$$\xymatrix{ 0\ar[r] & {\mathcal{O}_X^k}\ar[r] & E \ar[r] &
F\ar[r] & 0\\
0\ar[r] & {\mathcal{O}_X^{k^\prime}}\ar[r]\ar[u] &
\Psi_a^1({\widehat E}^\prime) \ar[r]\ar[u] &
\Psi_a^1({\widehat F}^\prime)\ar[r]\ar[u] & 0\\
& 0\ar[u] & 0\ar[u] & 0\ar[u] & }$$ with $\mu(\Psi_a^1(\widehat
F^\prime))=\mu(F)$. The stability of the original BGN extension
gives
\begin{equation}\label{bgn2}
\frac{k^\prime}{\hat r^\prime-a\hat d^\prime}-\frac{k}{r}>0
\end{equation}

Now we have \begin{align} \label{l1}&
r+ad=r(1+a\mu(E))=r\{1+a\,\left(\frac{r-k}{r}\right)\mu(F)\}\\
\label{l2} &\hat r^\prime-a\hat d^\prime=r(1-a\mu(\widehat
E^\prime))=\hat r^\prime\{1-a\,\left(\frac{\hat
r^\prime-k^\prime}{\hat r^\prime}\right) \mu(\widehat F^\prime)\}
\end{align} On the other hand the equality $\mu(\widehat
F^\prime)=\mu(\Phi_a^0(F))$ gives $\mu(\widehat
F^\prime)=\frac{\mu(F)}{1+a\mu(F)}$. If we substitute this in
$(\ref{l2})$, after some algebra, we get
\begin{equation}\label{l3}\hat r^\prime+a\hat d^\prime=\frac{\hat
r^\prime+a k^\prime\mu(F)}{\hat r^\prime(1+a\mu(F))}
\end{equation}

Substituting (\ref{l1}) in (\ref{bgn1}) and (\ref{l3}) in
(\ref{bgn2}) we get, after clearing denominators,

\begin{align}
&{k^\prime}(r+a(r-k)\mu(F))-{k}{\hat r^\prime}\leq 0\\
&{k^\prime{r}(1+a\mu(F))}-{k}{(\hat r^\prime+a k^\prime\mu(F))}>0
\end{align} Expanding this inequalities we get finally

\begin{align}
&{k^\prime}r-{k}{\hat r^\prime}+a{k^\prime}(r-k)\mu(F))\leq 0\\
& k^\prime{r}-{k}\hat r^\prime+ak^\prime({r}-{k})\mu(F)>0
\end{align} which can not hold simultaneously.

\end{proof}

Using similar arguments we prove:

\begin{thm}\label{large-psi}
The Fourier-Mukai transform $\Psi_a$ induces a map
$$\Psi_a^1\colon G_L(r+ad,d,k)\to G_L(r,d,k).$$
\end{thm}

As a consequence of Theorems \ref{large-fi}, \ref{large-psi}, and
the fact that $\Psi_a[1]$ and $\Phi_a$ are quasi-inverses we have:

\begin{thm}\label{lg-iso}
The Fourier-Mukai transform $\Phi_a$ induces an isomorphism of
moduli spaces
$$\Phi_a^0\colon G_L(r,d,k)\to G_L(r+ad,d,k)$$ Therefore, the
isomorphism type of $G_L(r,d,k)$  depends only on the class
$[r]\in\mathbb Z/d\,\mathbb Z$.
\end{thm}

\section{Birational type of the moduli spaces
$G(\alpha;r,d,k)$}\label{bir-type}

Let us denote by $0=\alpha_0<\alpha_1<\cdots<\alpha_L$ the
critical values for coherent systems of type $(r,d,k)$. As we have
mentioned above, the moduli spaces $G(\alpha;r,d,k)$ for any two
values of $\alpha\in (\alpha_i,\alpha_{i+1})$ coincide. Thus there
is only a finite number of different moduli spaces. On the other
hand one has (see \cite[Theorem 4.4]{LN}) the following important
result.

\begin{thm}\label{bir-t} The birational type of $G(\alpha;r,d,k)$ is
independent of $\alpha\in(\alpha_0,\alpha_L)$.
\end{thm}

Hence, the birational type shared by the finitely many different
moduli spaces can be determined by considering any of them. From
the point of view of the study carried out in previous sections by
means of Fourier-Mukai transforms, the natural candidates to
consider are $G_0(r,d,k)$ and $G_L(r,d,k)$. However, the latter is
in this respect less convenient than the former since in Section
\ref{lg-alpha} we have studied $G_L(r,d,k)$ under the restriction
$0<k<r$, whereas $G_0(r,d,k)$ has been considered in full
generality in Section \ref{sm-alpha}.

\begin{thm}
Let $a$ be a positive integer. The birational types of
$G(\alpha;r,d,k)$ and $G(\alpha;r+ad,d,k)$ are the same. Therefore,
the birational type of $G(\alpha;r,d,k)$ depends only on the class
$[r]\in\mathbb Z/d\,\mathbb Z$.
\end{thm}

\begin{proof}
By Theorem \ref{bir-t}, the birational types of $G(\alpha;r,d,k)$
and $G(\alpha;r+ad,d,k)$ are the same as those of $G_0(r,d,k)$ and
$G_0(r+ad,d,k)$, respectively. Due to Theorem \ref{sm-iso}, the
Fourier-Mukai transform $\Phi_a$ induces an isomorphism
$\Phi_a^0\colon G_0(r,d,k)\to G_0(r+ad,d,k)$. This finishes the
proof in view of Theorem \ref{bir-t}.
\end{proof}
\eject

\noindent \textbf{\Large Acknowledgements}

\noindent The authors would like to thank the Department of
Mathematical Sciences of the University of Liverpool for its
hospitality while some parts of this work were carried out, and
Professor P. E. Newstead for stimulating discussions.


\begin{thebibliography}{100}\frenchspacing

\bibitem{AHR} B. Andreas, D. Hern\'{a}ndez Ruip\'{e}rez, \emph{Fourier Mukai
transforms and applications to string theory}, Rev. R. Acad. Cien.
Serie A. Mat. vol. {\bf 99} (1), 2005, 29–77.

\bibitem{A} M. F. Atiyah, \emph{Vector bundles over an elliptic curve},
Proc. London Math. Soc. {\bf 7} (1957), 414-452.

\bibitem{BB} C. Bartocci, I. Biswas, \emph{Higgs bundles and the
Fourier-Mukai transform,} Southeast Asian Bull. Math. \textbf{25}
(2001),  201-207.

\bibitem{BBH2} C. Bartocci, U. Bruzzo, D. Hern\'andez Ruip\'erez, \emph{
A Fourier-Mukai transform for stable bundles on $K3$ surfaces}, J.
Reine Angew. Math. \textbf{486} (1997), 1--16.

\bibitem{BBHJ} C. Bartocci, U. Bruzzo, D. Hern\'{a}ndez Ruip\'{e}rez, M.
Jardim,  ``Fourier-Mukai and Nahm transforms in geometry and
mathematical physics'', book to appear in Progress in Mathematical
Physics, Birkh\"{a}user.

\bibitem{BBHM} C. Bartocci, U. Bruzzo, D. Hern\'{a}ndez Ruip\'{e}rez, J. M. Mu\~{n}oz
Porras, \emph{M. Mirror symmetry on $K3$ surfaces via
Fourier-Mukai transform,} Comm. Math. Phys. {\bf 195} (1998), no.
1, 79--93.

\bibitem{BO} A. Bondal, D. Orlov, \emph{Reconstruction of a variety from the
derived category and groups of autoequivalences}, Compositio Math.
\textbf{125} (2001), 327--344.

\bibitem{BG} S. B. Bradlow, O. Garc\'{\i}a-Prada, \emph{An application of
coherent systems to a Brill-Noether problem}, J. Reine Angew.
Math. \textbf{14} (2002), 123-143.

\bibitem{BGMN} S. B. Bradlow, O. Garc\'{\i}a-Prada, V. Mu\~noz and P. E.
Newstead, \emph{Coherent systems and Brill-Noether theory},
Internat. J. Math.  \textbf{14} (2003), 683-733.

\bibitem{BGN} L. Brambila-Paz, I. Grzegorczyk, P. E. Newstead, \emph{Geography
of Brill-Noether loci for small slopes,} J. Algebraic Geom.
\textbf{6} (1997),  645-669.

\bibitem{Br} T. Bridgeland, \emph{Fourier-Mukai transforms  for elliptic
surfaces,} J. Reine Angew. Math. \textbf{498} (1998), 115--133.

\bibitem{GHPT} O. Garc\'{\i}a-Prada, D. Hern\'{a}ndez Ruip\'{e}rez, F. Pioli, C. Tejero Prieto,
\emph{Fourier-Mukai and Nahm transforms for holomorphic triples on
elliptic curves,}  J. Geom. Phys.  \textbf{55}  (2005), 353-384.

\bibitem{HP}  D. Hern\'andez Ruip\'erez, and J.M.
Mu\~noz Porras, \emph{Stable sheaves on elliptic fibrations},
Journal of Geometry and Physics, {\bf 43} (2002) 163-183.

\bibitem{Hi} L. Hille, M. Van den Bergh, \emph{Fourier-Mukai transforms}, to
appear in the "Handbook on tilting theory".

\bibitem{KN} A. D. King, P. E. Newstead, \emph{Moduli of Brill-Noether pairs on algebraic
curves,} Internat. J. Math. \textbf{6} (1995), 733-748.

\bibitem{LN} H. Lange, P. E. Newstead, \emph{Coherent
systems on elliptic curves,}  Internat. J. Math. \textbf{16}
(2005), 787-805.

\bibitem{LP} J. Le Potier, \emph{Syst\`emes coh\'erents et structures de
niveau,} Ast\'erisque \textbf{214} (1993).


\bibitem{Mac} A. Maciocia, \emph{Gieseker stability and the Fourier-Mukai transform for
abelian surfaces,}  Quart. J. Math. Oxford Ser. (2)  \textbf{47}
(1996), no. 185, 87--100.

\bibitem{Mu} S. Mukai, \emph{Duality between $D(X)$ and $D(\hat{X})$ with
its application to Picard sheaves,} Nagoya Math. J. {\bf 81} (1981), 153--175.


\bibitem{Mu1}S. Mukai, \emph{Fourier functor and its application to the moduli
of bundles on an abelian variety}, Adv. Pure Math. {\bf 10} (1987),
515-550.

\bibitem{O1} D. O. Orlov, \emph{Equivalences of derived
categories and $K3$ surfaces,} J. Math. Sci. (New York)
\textbf{84} (1997), no. 5, 1361-1381.

\bibitem{O2} D. O. Orlov, \emph{Derived categories of coherent sheaves on abelian
varieties and equivalences between them,} Izv. Ross. Akad. Nauk
Ser. Mat. \textbf{66} (2002), no. 3, 131-158.

\bibitem{RV} N. Raghavendra, P.A. Vishwanath, \emph{Moduli of pairs and
generalized theta divisors}, T\^ohoku Math. J. \textbf{46} (1994),
321-340.

\bibitem{Tu} L. W. Tu, \emph{Semistable bundles over an elliptic
curve}, Adv. in Math.{\bf 98} (1993) 1-26.

\bibitem{Yo} K. Yoshioka, \emph{Moduli spaces of stable sheaves on
abelian surfaces}  Math. Ann.  321  (2001),  no. 4, 817--884.

\end{thebibliography}
\end{document}